\documentclass[a4paper,11pt, reqno]{amsart}

\usepackage{amsmath,amsfonts,amsbsy,amsgen,amscd,mathrsfs,amssymb,amsthm}

\newtheorem{definition}{Definition}[section]
\newtheorem{theorem}{Theorem}[section]
\newtheorem{lemma}[theorem]{Lemma}
\newtheorem{corollary}[theorem]{Corollary}

\newtheorem{remark}{Remark}[section]

\providecommand{\abs}[1]{\lvert#1\rvert}
\providecommand{\norm}[1]{\lVert#1\rVert}

\begin{document}

\title[Exponential Attractors in Hilbert and Banach Spaces]
{Rotating Navier-Stokes-$\alpha$ equations: Exponential Attractors in Hilbert and Banach Spaces}

\author{Bong-Sik Kim}
\address{P.O. Box 10021, Department of Mathematics \& Natural Sciences,
 American University of Ras Al Khaimah,  
United Arab Emirates}
\email{bkim@aurak.ac.ae}


\begin{abstract}
  This article covers the construction of exponential attractors  
  in two different functional space settings; one is in Hilbert's
space, and the other is in the Banach space. The former relies on the
squeezing properties of solution trajectories, but the latter does not. 
We present these different methods for constructing exponential 
attractors using the three-dimensional rotating
 Navier-Stokes-$\alpha$ equations. 
\vspace{2mm}


\vspace{2mm}


\end{abstract}



\maketitle

\section{Introduction}

In attempting to model phenomena in nature that change with time, 
the model equations generally come in as a system of partial differential equations 
and nonlinearities occur in the modeling process.  The obtained nonlinear system evolves in time 
and exhibit gradual or rapid change as time proceeds. Typically, it has infinite dimensional 
aspects, the dimensions here being the number of parameters which is necessary 
to describe the configuration of the motion 
at a given instant in time.

If there is no restoring force, the flow of a quantity, such as density, concentration, or heat, 
 tends to dissipate. 
 The dissipative effect is reflected in an infinite-dimensional nonlinear system, 
 and defines a forward regularizing flow in an adequate phase space $X$,
 $S(t):X\rightarrow X$,
  containing an \textit{absorbing set}.
The absorbing set $B\subset X$ is a bounded set that attracts all bounded solutions in
\textit{finite time}. The existence of such an absorbing set can be taken as a definition of
dissipative partial differential equations.

Since all solution trajectories of dissipative systems  eventually enter and stay
in $B$, we may expect the existence of a set which would capture all the asymptotic dynamics.
Such a set is called the \textit{global attractor} and it is the largest set that enjoys 
positively and negatively invariant properties under the flow. More precisely,
the global attractor $\mathcal{A}\subset X$ is 
\begin{itemize}
\item[(i)] the maximal compact invariant set, 
 $S(t)\mathcal{A} = \mathcal{A}\ \mbox{for all}\ t\geq 0$; and
\item[(ii)] the minimal set that attracts all bounded sets,
$ \mathrm{dist}(S(t)B, \mathcal{A}) \rightarrow 0
  \ \mbox{as}\ t\rightarrow \infty $
 for any bounded set $B\subset X$.
\end{itemize}
If a dynamical system posesses a global attractor, it is unique for the system.
The global attractor, however, is not stable under pertubations of
the underlying evolution equations and the attraction rate
can be arbitrarily  slow. Those reasons led to the development of
the concept of exponential attractors, 
first introduced in \cite{Eden} in the Hilbert space context, 
and further generalized to Banach spaces in \cite{Dung}.
In contrast to the global attractor, exponential attractors are strongly stable, 
attract all solution trajectories at exponential rates, but not unique.
Further, 
the attractors often have finite fractal dimensions and the
asymptotic behavior of the given system
can be approximated by a finite-dimensional dynamical system.
In its numerical interpretation,
 the existence of a finite-dimensional attractor 
 guarantees that long-time behavior of the given system
 can be numerically approximated using a discrete system with
 a finite number of degrees of freedom.

We can  configure exponential attractors for a dissipative dynamical system 
in two different functional settings. 
One is to construct the attractors in the Hilbert space setting, 
and the other is to build th
  Here, of course, it is the lack of conservation of classical energy for
  the unfiltered $V_{\alpha}$ that creates the difficulty. 
Based on the existence and regularity results we
establish the existence of exponential attractors of the 3D
RNS-$\alpha$ equations.e attractors in the Banach space setting. 
The former implicitly relies on some squeezing properties of trajectories \cite{Eden}. 
No squeezing conditions are needed in the latter. We demonstrate the two different ways by 
constructing the exponential attractors of the three-dimensional 
Lagrangian-averaged Navier-Stokes equations for uniformly rotating fluid flows.

 The Navier-Stokes-$\alpha$
equations (also known as the Lagrangian-averaged Navier-Stokes equations)
  were introduced as a turbulence 
closure model
in 1998, \cite{Holm1}. 
  Here, of course, it is the lack of conservation of classical energy for
  the unfiltered $V_{\alpha}$ that creates the difficulty. 
Based on the existence and regularity results we
establish the existence of exponential attractors of the 3D
RNS-$\alpha$ equations.
This work is based on theoretical results
from \cite{Kim1, Kim2}, where the Navier-Stokes-$\alpha$ equations are
considered for fluids in a periodic box, with uniform rotation about the vertical
axis $e_3=(0, 0, 1)$ of angular frequency $f=2\Omega$. In a rotating frame of reference,
the Rotating Navier-Stokes-$\alpha$ equations (RNS-$\alpha$ equations) 
are given by
\begin{eqnarray}
  \frac{\partial V}{\partial t}+
        (U\cdot\nabla)V
         + V_j\nabla U^j
         +fe_3\times U
     & = & -\nabla \pi+\nu\Delta V+F\label{RLANS1} \\
             \nabla\cdot V& = &  0 \label{RLANS2}\\
         V(t,x)|_{t=0}&=&
         V(0,x)=V(0),\label{RLANS3}
 \end{eqnarray}
where
\begin{center}
   \begin{eqnarray}
     V(t,x) &=& (V_1,V_2,V_3)\ \ \ \mbox{the velocity vector}, \label{RLANS4} \\
     U(t,x) &=& (I-\alpha^2\Delta )^{-1}V(t,x) \ \ \ \mbox{the filtered
       velocity}, \label{RLANS5} \\
      \pi    &=&
        \frac{p}{\rho}-\frac{1}{2}|U|^2
        -\frac{\alpha^2}{2}|\nabla U|^2  \ \ \mbox{the modified pressure}.\label{RLANS6}
   \end{eqnarray}
\end{center}
Here $x=(x_1,x_2,x_3)$, $f=2\Omega$ is the Coriolis parameter,
$F=(F_1,F_2,F_3)$ is a divergence free force, $\nu >0$ is
the kinetic viscosity, $\rho$ is the fluid density, and $p$ is the pressure.
For simplicity we will assume the forcing term to be time independent; that is,
$F(x,t)\equiv F(x)$. the parameter $\alpha$ is a length scale,  
below which wave acitivity is filtered, with
$0<\alpha\ll 1$. 

Main difficulty for this 3D flow model is in maintaining
all the estimates bounded when $\alpha\rightarrow 0^+$ during the construction of 
the exponential attractors. 
Ilyin and Titi \cite{Ilyin1} estimated attractor dimensions for two-dimensional Navier-Stokes-$\alpha$
  equations. Their estimates, however, blow up as $\alpha\rightarrow 0^+$. 
  Gibbon and Holm \cite{Gibbon1} obtained length-scale estimates for NS-$\alpha$ equations in 
  terms of the Reynolds number, which blow up in the limit when $\alpha\rightarrow 0^+$, too. 
  Several other time-averaged estimates related 
  to NS-$\alpha$ equations don't remain finite in the limit
   (Table 1 in \cite{Gibbon1}). 
  They analyzed the system in the context of the filtered velocity $U = (I-\alpha^2\Delta )V$.
  Instead,  we study the system from the perspective of 
   the non-filtered velocity $V = (I-\alpha^2\Delta )^{-1}U$.
  The Helmholtz inverse operator $\mathcal{R}_\alpha = (I-\alpha^2\Delta )^{-1}$ plays a  crucial role
   in the process, leading to 
  uniform estimates that remain finite  as $\alpha\rightarrow 0^+$.
  
Another difficulty we encounter is that
the 3D RNS-$\alpha$ equations lack  a uniform spatial $L^2$-norm of 
the unfiltered
velocity $V$.  The best uniform in $\alpha$ estimates
is restricted to $H^{\beta}, \beta > 5/2$. 
This leads to the construction of 
an absorbing ball in a weaker topology than the topology
    of the initial data set, which is a typical feature of the infinite-dimensional dynamical
  systems methods applied to damped hyperbolic PDE's (Ch6, \cite{Eden}).
   This notion enables
  the construction of exponential attractors in the
  ``$H^{\beta}$-$H^{\gamma}$" sense.

Let's start with the definition of exponential attractors(\cite{Eden}):
\begin{definition}\textbf{(Exponential Attractor)}
Let $(E, d)$ be a complete metric space with a metric $d$, $X$ a compact subset of $E$, and
$\{ S(t) | t\geq 0\}$ the semigroup on $X$ for the topology of $E$.
Assume that $S(t)$ possesses a global attractor $\mathcal{A}$.
 A compact set $\mathcal{M}$
 is called an exponential attractor for the semidynamical system  $(S(t),X)$ if
 \begin{itemize}
    \item[(i)\,] $\mathcal{A}\subseteq \mathcal{M}\subseteq X$,
   \item[(ii)] $S(t)\mathcal{M}\subseteq\mathcal{M}$ for $t\geq 0$,
           (positively invariant under the flow),
   \item[(iii)] the fractal dimension of $\mathcal{M}$ is finite,
              $\dim_F(\mathcal{M})<\infty$, and
   \item[(iv)]there exists positive
             constants $c_0$ and $c_1$ such that
         \[ d_h(S(t)X, \mathcal{M}) \leq c_0e^{-c_1t},\ \forall t\geq 0, \]
         where $d_h$ is the Hausdorff semi-distance for the metric $E$ defined by
        \[ d_h(X,Y)=\sup_{x\in X} \inf_{y\in Y} d(x,y). \]
 \end{itemize}
\end{definition}

 \section{Existence of an absorbing set}

We denote $P_L$ as the usual Leray projector and
introduce an operator
$\mathcal{R}_{\alpha}=(1-\alpha^2\Delta )^{-1}$, which is
defined by $\mathcal{R}_{\alpha}v=(1-\alpha^2\Delta )^{-1}v$.
We also define a  bilinear operator $B_{\alpha}$ on divergence free vector fields by
\begin{equation}\label{bilinear}
B_{\alpha}(u,v)=P_L[ ( \mathcal{R}_{\alpha}u\cdot\nabla ) v
+v_j\nabla (\mathcal{R}_{\alpha}u)_j ].
\end{equation}
Then (\ref{RLANS1}) takes  the form
\begin{equation}\label{RLANS.Leray}
  \frac{\partial V}{\partial t}+fP_LJP_L\mathcal{R}_{\alpha}V+\nu AV+
    B_{\alpha}(V,V) = F,
\end{equation}
where $A=-P_L\Delta$ is the Stokes operator and $J$ is a rotation matrix given by
\[ J = \left( \begin{array}{ccc}
                  0 & -1 & 0 \\
          1 & 0 & 0 \\
          0 & 0 & 0
        \end{array} \right)
 \]
The system is considered subject to periodic boundary conditions
in a lattice $Q=[0, 2\pi a_1]\times[0, 2\pi a_2]\times [0, 2\pi a_3]$
as well as stress-free boundary conditions in the vertical.
The corresponding function spaces are  Fourier-Sobolev spaces of periodic functions,
$H^s$, $s\geq 0$, with the norm
\begin{equation*}
  \norm{u}_s^2=\sum_{n\in Z^3}\abs{\check{n}}^{2s}\abs{u_n}^2,
\end{equation*}
where  $n=(n_1,n_2,n_3)\in Z^3$ is a wave number and
 $\check{n}=(\check{n}_1,\check{n}_2,\check{n}_3)$ with $
 \check{n}_j=n_j/a_j$ for $j=1,2,3$.
We set $a_1=1$ without loss of generality.

The existence of unique regular solutions for all $f = 2\Omega$ greater than some threshold $f^*$
has been proved in \cite{Kim1}, which has led to the existence of absorbing sets:

\begin{theorem}\label{theorem22} (Theorem 6.7 in \cite{Kim1})
    Let $0\leq \alpha ,\ \nu >0; $ let $a_1, a_2, $ and $a_3$ be arbitrary and fixed.
 Let $\beta > 5/2, \gamma > \beta +4,$ and
 $F$  a  time-independent force such that
\begin{equation*}
  \norm{F}_{\beta -1}^2 \leq M_{\beta F}^2  \ \ \mbox{and} \ \
  \norm{F}_{\gamma -1}^2 \leq M_{\gamma F}^2.
\end{equation*}
Let $V(0)\in B_{\gamma I}$ be  initial data in a ball in $H^{\gamma}$.
Let $diam (B_{\gamma I}) < 2\rho_{\gamma I}$ in $H^{\gamma}$ norm and
{\small $diam (B_{\gamma I}) < 2 \rho_{\beta I}$} in $H^{\beta}$ norm.
   Then, for each {\small $f\geq f^{*}(M_{\beta F}, M_{\gamma F},
   \rho_{\gamma I}, \rho_{\beta I}, \nu , a_1, a_2, a_3)$},
  the three-dimensional (3D) rotating Navier-Stokes-$\alpha$ equations 
  possess an absorbing set $B_{\beta}$
  in $H^{\beta}$; that is, there exists $t_{\beta}=t_{\beta}(\rho_{\beta I})$, such that
$f\geq f^{*}$ and $V(0)\in B_{\gamma I}$ imply
  \begin{equation*}
    \norm{V(t)}_{\beta} \leq \rho_{\beta} \ \ \ \mbox{and}\ \ \
    \nu \int_{t}^{t+1}\norm{V(\tau )}_{\beta +1}^2\,d\tau \leq M_{\beta+1}^2,
  \end{equation*}
for all $t\geq t_{\beta}$.
   This absorbing set is uniform in $\alpha$ and $\rho_{\beta}=\rho_{\beta}(M_{\beta F}, \nu , a_1, a_2, a_3),$
  $ M_{\beta +1}=M_{\beta +1}(M_{\beta F}, \nu , a_1, a_2, a_3)$ (with no dependence on
   $M_{\gamma F}$,
   nor on $\rho_{\gamma I}$).
\end{theorem}

As $\beta > 5/2$, the semiflow $S_{\alpha}(t)$ is compact on $H^\beta$ and
we can take $B_\beta$ compact in $H^\beta$, modulo a small translate in time.

\begin{remark}
 Existence of unique regular solutions of the exact rotating Navier-Stokes equations ($\alpha =0$)
 was established by Babin, Mahalov and Nicolaenko in \cite{BMN1} and \cite{BMN2}.
\end{remark}

 \section{Existence of exponential attractors in Hilbert spaces}

The existence of exponential attractors for the system (\ref{RLANS1})-(\ref{RLANS3})
in Hilbert spaces
was established in \cite{Kim1}. The procedure and results
are reproduced in this section. We first would like to point out that exponential attractors are, 
unlike a global attractor, stable under perturbations of the underlying evolution equations.
The full 3D rotating Navier-Stokes systems, including Lagrangian-averaged Navier-Stokes-$\alpha$
equations, are considered to be an $f$-singular perturbation from $f$-singular limit equations.
With manifolds that stay stable under the perturbation, we are able to talk about
convergence as $f\rightarrow \infty$. See \cite{Kim1} for more details on this. Now,
we start out by recalling the procedure with which exponential attractors are 
constructed in Hilbert spaces. 

Let $E$ be a Hilbert space with norm $\norm{\cdot}_E$ induced by the
inner product $(\cdot ,\cdot )_E$. Let $X$ be a compact subset of $E$
and $S:X\rightarrow X$ a Lipschitz  continuous map with Lipschitz constant $L$.
Then $S$ possesses a global attractor $\mathcal{A}$
which is a compact, connected set given by
\begin{equation*}
  \mathcal{A}=\bigcap_{n=1}^{\infty}S^n(X)
\end{equation*}
(Theorem 2.4.2, \cite{Hale}). 
Exponential attractors for a map $S$ are defined as
\vspace{.1in}
\begin{definition} \textbf{(Discrete Exponential Attractor) }
  A compact set $\mathcal{M}$ is called an exponential attractor for
  $(S, X)$ if $\mathcal{A}\subset \mathcal{M} \subset X$ and
  \begin{itemize}
    \item[(i)\,] (positively invariant) $S(\mathcal{M})\subset \mathcal{M}$.
    \item[(ii)] $\mathcal{M}$ has finite fractal dimension,
                  $dim_F(\mathcal{M})<\infty$.
    \item[(iii)] There exist positive constants $c_0$ and $c_1$ such that
        \begin{equation}
          d_h(S^nX,\mathcal{M})\leq c_0e^{-c_1n},\ \ \forall n\geq 1 \nonumber
        \end{equation}
	where $d_h$ is the standard Hausdorff semi-distance between two sets.
  \end{itemize}
\end{definition}
\vspace{.1in}

In establishing the existence of discrete exponential attractors
key techniques are those based on examining the difference of two solutions and
verifying the squeezing property on the underlying mapping $S$. 
The idea of the squeezing property is that we can split the phase space $X$
into a finite-dimensional subspace and its infinite-dimensional
orthogonal complement, such that the finite-dimensional part of the solution dominates; or
if not, then at least the solutions are closer together than they were
at $t=0$, which serves to dampen the effect of
such \textit{ill-behaved} solutions:

\vspace{.1in}
\begin{definition}
  Let $E$ be a Hilbert space and $X$ a subset of $E$. 
  A map $S$ has 
  the squeezing property in $X$ if, for some $\delta\in (0, \frac{1}{4})$,
  there exists an orthogonal projection $P_{N_0}=P_{N_0}(\delta )$ of finite rank $N_0=N_0(\delta )$
  such that, $\forall u, v\in X$,  if
  $ || (I-P_{N_0})(Su-Sv)||_E \geq ||P_{N_0}(Su-Sv)||_E$ then
  $||Su-Sv||_E \leq \delta ||u-v||_E$.
\end{definition}
\vspace{.1in}

In general, to decrease $\delta$ we need to increase the rank of the orthogonal projection
$P_{N_0}$ (that is, the dimension of $P_{N_0}X$).
The squeezing property guarantees the existence 
of discrete exponential attractors (Ch. 2, \cite{Eden}):
\begin{theorem}\label{thm7.3}
 If $S$ has the squeezing property in $X$, then there exists an exponential attractor
 $\mathcal{M}$ for $(S,X)$ and, moreover,
 \begin{equation}
   d_B(\mathcal{M})\leq N_0 \, 
   \max \left\{ 1, \log (\frac{2L}{\delta}+1)/\log (\frac{1}{\theta})\right\}, \nonumber
 \end{equation}
 where $\theta\in (4\delta ,1)$ arbitrary and $d_B$  is the fractal box dimension for
 the metric $E$.
\end{theorem}
\vspace{.1in}

We now turn to the continuous case.
 Given the semigroup $\{ S(t) | t\geq 0\}$ of solution operators, we will choose a positive $t_*$ small enough such that
$S_*=S(t_*)$ possesses the squeezing property in $X$.
 If $S_*$ is Lipschitz continuous, then the existence of a discrete exponential attractor
$\mathcal{M}_*$ for $(S_*,X)$ is guaranteed  by Theorem \ref{thm7.3}.
 Next we define
\[ \mathcal{M} = \bigcup_{0\leq t\leq t_*}S(t)\mathcal{M}_* \]
and $G: [0, T]\times\mathcal{M}_*\rightarrow\mathcal{M}$ as $G(t,x)=S(t)x$.
If $G$ is Lipschitz, then it can be shown that $\mathcal{M}$ is a compact set with finite
fractal box dimension, and $\mathcal{M}$ will be an exponential attractor for
$(S(t),X)$ (Theorem 3.1, \cite{Eden}).
The exponential attractors for the continuous dynamical systems generated by 
a semigroup $\{ S(t)\}_{t\geq 0}$ are unions of exponential attractors restricted by
squeezing time $t_*$. In addition, given an estimate for $\mathcal{M}_{*}$,
it is not difficult to get an estimate for the fractal box dimension of $\mathcal{M}$
 (Theorem 3.1, \cite{Eden}):
\vspace{.1in}
\begin{theorem}
  Let $S(t_*)$ have the squeezing property in $X$ for some time $t_*>0$ and let
  $\mathcal{M}$ be an exponential attractor for $(S(t),X)$ and $G(t_*,x)=S(t_*)x$
  for $x\in X,\ t\geq 0$. If $G(t_*,\cdot )$ is Lipschitz in $X$ with Lipschitz
  constant $L_*$, 
  then 
  \[ d_B(\mathcal{M})\leq d_B(\mathcal{M}_*)+1. \]
  Furthermore, 
  \[ d_h(S(t)X, \mathcal{M}) \leq cL_*\exp\left( \frac{-(\ln 8) t}{t_*}\right) \]
  for all $t\geq 0$, where $c$ is a positive constant.
\end{theorem}
\vspace{.1in}

Now we follow the above procedure to establish the existence of an 
exponential attractor in $L^2$ for the 3D RNS-$\alpha$ equations.
We do this for all $f$ that allow the existence of a global attractor.
Assume that $F$ is time-independent and smooth and that
$f\geq f^*$ as in Th \ref{theorem22}. Let $S_{\alpha}(t)$ be the semiflow for solutions of
the 3D RNS $\alpha$-equations and let $B_{\beta}, \beta > 5/2$, be
the compact absorbing set obtained in Theorem \ref{theorem22}. Set
\begin{equation}\label{sec7.1}
 X_{\alpha ,\beta} =
 \overline{\cup_{t\geq t_{\beta}(B_\beta )+\frac{1}{\nu \lambda_1}}
           S_{\alpha}(t)B_{\beta}}^{|\cdot |}, \nonumber
\end{equation}
where the closure is taken in $L^2$-topology and $\lambda_1$ denotes
the first eigenvalue of the Stokes operator. Then $X_{\alpha ,\beta}$
is a bounded subset of $B_{\beta}$, compact in $H^s, 0\leq s<\beta$, and positively invariant
under $S_{\alpha}(t)$ such that, for all $V_{\alpha}(0)\in X_{\alpha
,\beta}$,
\begin{equation}\label{sec7.2}
 ||S_{\alpha}(t)V_{\alpha}(0)||_{H^\beta} \leq \rho_{\alpha ,\beta},
  \ \ \forall t\geq 0, \nonumber
\end{equation}
where $\rho_{\alpha ,\beta}$ 
is the uniform bound obtained in Th \ref{theorem22}.
In particular, there exist absolute bounds
 $\rho_{\alpha ,s}=\rho_{\alpha ,s}(M_{\beta F}, \nu , a_1, a_2, a_3)$
such that $\norm{S_{\alpha}(t)V_{\alpha}(0)}_{H^s}\leq\rho_{\alpha ,s}
 \leq \rho_{\alpha ,\beta}$ for $0\leq s<\beta$.
We will denote $\rho_H$ for  $\rho_{\alpha ,0}$ and $\rho_V$ for $\rho_{\alpha ,1}$.
Since $X_{\alpha ,\beta }$ is compact in $H^s$ for $0\leq s<\beta $, we can
deduce that the underlying semigroup $S_{\alpha}(t)$ is \textit{uniformly compact} for large $t$
 so that 
it possesses a unique global attractor $\mathcal{A}$ in $H^s$ for $0\leq s<\beta $,
(Theorem 1.1, \cite{Temam1}). Moreover, it can be proved that
$\mathcal{A}$ lies in $H^{\beta}$ for $\beta > 5/2$.

We consider the solution operator $S_{\alpha}(t)$ as a map from $X_{\alpha ,\beta}$
into $X_{\alpha ,\beta}$. We only need to show that there exists a squeezing time
$t_*$ such that the discrete operator $S_{*}=S_\alpha (t_*)$ has the
squeezing property in $L^2$-topology. To achieve it we first examine
the difference between two solutions, $V_a$ and
$V_b$, of 3D RNS-$\alpha$
 equations in $X_{\alpha ,\beta}$.
Let $W=V_a-V_b$ and $W^\prime = \frac{V_a+V_b}{2}$. Then
$W$ satisfies the equation
\begin{eqnarray}
  \frac{\partial W}{\partial t}+\nu AW+fM\mathcal{R}_{\alpha}W =
     -\left[ B_{\alpha} (W^\prime ,W)+B_{\alpha}(W,W^\prime )\right] && \label{eq7.3}  \\
 W(0) = V_a(0)-V_b(0). \hspace{1.5in} &&
\end{eqnarray}
Taking the inner product with $2W$ yields
\begin{equation}\label{eq7.5}
  \frac{d}{dt}\abs{W}^2+2\nu\norm{W}^2 \leq 2\left\{ \, \abs{<B_{\alpha}(W^\prime ,W), W>}
     + \abs{<B_{\alpha}(W ,W^\prime ), W>} \right\} \, ,
\end{equation}
where $B_{\alpha}(u ,v)=(\mathcal{R}_{\alpha}u\cdot\nabla )v+v_j\nabla
  (\mathcal{R}_{\alpha}u)_j$.
Estimating the right hand side of (\ref{eq7.5}) and using  Young's inequality yield
\begin{equation}\label{eq7.6}
  \frac{d}{dt}\abs{W}^2+\nu\norm{W}^2 \leq
     \frac{K_1}{\nu^3}\abs{W}^{2}, \hspace{.3in}
\end{equation}
where $K_1=c_1^4\rho_V^4$ with $c_1$ a constant. Letting $\lambda
(t)=\frac{\norm{W(t)}^2}{\abs{W(t)}^2}$, (\ref{eq7.6}) becomes
\begin{equation}
  \frac{d}{dt}\left[ \ln \abs{W(t)}^2\right] \leq
        -\nu\lambda (t)+\frac{K_1}{\nu^3} \nonumber
\end{equation}
so that
\begin{equation}\label{eq7.7}    
    |W(t)|^2 \leq \delta (t)|W(0)|^2
\end{equation}
with
\[
  \delta (t)=\exp \left(
        -\nu \int_0^t\lambda (s)\,ds +\frac{K_1}{\nu^3}t\right).
\]

Next, we need to find a time $t_*$ such that the estimate for $\delta (t_*)$ allows
squeezing. Thus it is essential to bound $\int_0^{t_*}\lambda (s)\,ds$, and 
following the exact line of section 6.1 in \cite{Markus} we obtain
\begin{equation}\label{eq7.10}
   t_* = \frac{c_3^2}{c_2}\frac{\nu^{3/2}}{K_2 K_3},
\end{equation}
where $K_2=c_2\rho_V$ and $K_3^2=\frac{27c_3^4}{2\nu^3}\rho_V^6+\frac{2}{\nu\lambda_1}\rho_H$
with $c_2$ and $c_3$  constants.
Furthermore,
\begin{equation}
  \int_0^{t_*}\lambda (t)dt \geq
     c_4\lambda_{N_0+1}\frac{\nu^{3/2}}{K_2K_3}, \nonumber
\end{equation}
where $c_4=\frac{1}{2}[1-\exp (-c_3^2/c_2)]>0$, so that 
\begin{equation}\label{eq7.11}
  \delta (t_*) \leq \exp \left( -\frac{c_4}{c_2}\lambda_{N_0+1}
   \frac{\nu^{5/2}}{K_3\rho_V}+
   \frac{c_5\rho_V^3}{\nu^{3/2}K_3}
   \right) ,
\end{equation}
where $c_5=\frac{27}{16}{c_1^4c_3^2}{c_2^2}$.
 By the definition of $K_3$ there exists a constant $\tilde{c}>0$
 such that
 \[ K_3 \leq \tilde{c}\left(
 \frac{\rho_V^3}{\nu^{3/2}}+\nu^{1/2}\lambda_1^{1/2}\rho_H\right)
 .\]
 Choosing $N_0$ such that
 \begin{equation}
   N_0 \geq \tilde{c}^{3/2}
    \max \left\{ \frac{1}{\lambda_1^{3/4}}
      \frac{(\rho_H\rho_V)^{3/2}}{\nu^3},
       \frac{\rho_V^6}{\lambda_1^{3/2}\nu^{6}}\right\} , \nonumber
 \end{equation}
 gives $\delta (t_*)<\frac{1}{8}$.
 Under the above condition of $N_0$, the following Lemma assures
 the existence of an exponential
 attractor $\mathcal{M}_{0}^*$ for $(S_*,X_{\alpha ,\beta})$  for
 $f\geq f_*$ (Ch 3, \cite{Eden}; Proposition 2.2.7, \cite{Markus}):
 \vspace{.1in}
\begin{lemma}\label{lemma7.7.2}
  Let $t_*>0$ be given and $u, v \in X$. Define
  \[ \lambda_* = \frac{\norm{w_*}^2}{\abs{w_*}^2} , \]
  where $w_*=S_*u-S_*v$. Then $S_*$ possesses the squeezing property in $X$, if there exists
  $\delta\in (0, 1/4)$ and $N_0=N_0(\delta )\in \mathcal{N}$,
  such that $\lambda_*>\frac{1}{2}\lambda_{N_0+1}$ implies that
  $\abs{S_*u-S_*v}<\delta\abs{u-v}$, for all $u, v \in X$.
\end{lemma}
\vspace{.1in}
 
 Furthermore, the Lipschitz constant for $S_*$ on $X_{\alpha ,\beta}$
 is estimated  as
 \begin{equation}
  L_*=\delta (t_*) \leq \exp \left(
     \frac{c_{5}\rho_V^3}{\nu^{3/2}K_3}\right) ,\nonumber
 \end{equation}
and hence
\begin{eqnarray*}
  d_h\left( S_{\alpha}(t)X_{\alpha ,\beta},\mathcal{M}_{0}^*\right)
    &\leq& cL_*\left( (\delta (t_*))^{1/t_*}\right)^t \\
    &\leq& cL_*\left( e^{-\ln 8}\right)^{t/t_*} \\
    &=& c_{\alpha F} e^{-\delta_{\alpha F}t},
\end{eqnarray*}
where $c_{\alpha F}=cL_*$ and $\delta_{\alpha F}=\frac{\ln 8}{t_*}$.

Now we summarize the results:
\vspace{.1in}
\begin{theorem}\label{theorem7.1}
  Let $F$ be a smooth, time-independent force and let $a=(a_1,a_2,a_3)$
  be a domain size parameter. For $f\geq f_*$ as in Th \ref{theorem22},
  let $X_{\alpha ,\beta}$ be the positively invariant set from
  (\ref{sec7.1}). Then $\{ S_{\alpha}(t) | t\geq 0\}$ restricted
  to $X_{\alpha ,\beta}$ admits an exponential attractor $\mathcal{M}_{0}$
  in $L^2$.
  Moreover, the rate of convergence to the exponential attractor is given by
  \[ d_h( S_{\alpha}(t)X_{\alpha ,\beta}, \mathcal{M}_{0}) \leq
    c_{\alpha F}e^{-\delta_{\alpha F}t},\]
  where $c_{\alpha F}, \delta_{\alpha F}$
  are constants, which only depend on $\nu, a, \rho_{\alpha H},\rho_{\alpha V}$
  and are independent of the angular frequency $f\geq f_0$ and $\alpha > 0$.
\end{theorem}
\vspace{.1in}

\begin{remark} $\mathcal{M}_{0}$ is bounded in $H^\beta$ and attracts all orbits
in the $L^2$-norm topology. It is compact in the space $H^\gamma , 0\leq \gamma < \beta$.
\end{remark}

\section{Existence of exponential attractors in Banach spaces}

Since the Hilbert space is also a Banach space, 
we may construct an exponential attractor using
the method developed by Le Dung and Nicolaenko in \cite{Dung}, which doesn't require 
the squeezing propertis of trajectories.
Let $\mathcal{L}(E)$ be the space of bounded linear maps from $E$ into itself.
For a given positive real $\lambda$ we denote by $\mathcal{L}_{\lambda}(E)$ the set of maps
$L\in\mathcal{L}(E)$ such that $L$ can be decomposed as $L=K+C$ with $K$ compact
and $\norm{C}<\lambda$. Here $\norm{C}$ denotes the norm of the operator $C$.
The following theorems were established in \cite{Dung}.
\begin{theorem}\label{theorem1}
  If there exists $\lambda\in (0,1)$ such that $D_xS(x)\in\mathcal{L}_{\lambda}(E)$ for
  all $x\in X$ then the discrete dynamical system $\{ S^n\}_{n=1}^{\infty}$
  possesses an exponential attractor.
\end{theorem}

Once the existence of exponential attractors for the discrete case is proved the result
for the continuous case follows in a standard way ( see Ch. 3, \cite{Eden}).
Define $S_{*}$ as the map induced by Poincar\'{e} sections of a Lipschitz continuous semiflow
$S(t),\ t\geq 0$ at the time $t=t^*$ for some $t^*>0$; that is, $S_{*}:=S(t^*)$.
Let $\{ S_{*}^n\}_{n\geq 0}$ be the discrete semigroup generated by $S_{*}$. Then
\begin{theorem}\label{theorem2}
  Let $X$ be a compact absorbing set for a continuous semiflow $S(t)$. Suppose that there is
  $t^*>0$ such that $S_{*}=S(t^*)$ satisfies the condition of Theorem \ref{theorem1}.
  Assume further that the map $G(x,t)=S(t)x$ is Lipschitz from $[0,T]\times X$ into $X$ for any
  $T>0$. Then the flow $\{ S(t)\}_{t\geq 0}$ admits an exponential attractor
  $\mathcal{M}$ as well as a unique global attractor $\mathcal{A}$.
\end{theorem}

Theorems  \ref{theorem1} and \ref{theorem2} were already proved by Temam (\cite{Temam1})
 and J. Hale (\cite{Hale}) for the
global attractor. 

\begin{theorem}\label{theorem27}
  Let $F$ be a smooth, time-independent force.
  The 3D RNS-$\alpha$ equations possess for $f>f^{*}$, where $f^{*}$ is defined in
  Theorem \ref{theorem22}, a global compact attractor $\mathcal{A}_{\beta}$ in
  the topology of $H^{\beta},\ \beta > 5/2$, as well as  exponential
  attractors $\mathcal{M}_{\beta}$ in the absorbing set $B_{\beta}$ established
  in Theorem \ref{theorem22}. Both fractal dimensions and rates of exponential
  attraction are uniform in $\alpha$.
\end{theorem}
\noindent\textbf{Proof.}
  We place ourselves in the context of the compact absorbing ball $B_{\beta}$ of
  Theorem \ref{theorem22}, which is absorbing in the $H^{\beta}$-topology
  the initial set $B_{\gamma I}$ in $H^{\gamma}$. We prove in the below Lemma
  \ref{lemma28}
  that the map $F(v,t)=S(t)v$ is Lipschitz from $[0,\ T]\times B_{\beta}$ into
  $B_{\beta}$ for any $T>0$, as well as the uniform Fr\'{e}chet Differentiablility
  of $S_{\alpha}(t)v$ with respect to $v\in B_{\beta},\ 0\leq t\leq T$.
  Then, all assumptions in Theorems \ref{theorem1} and \ref{theorem2} are
satisfied uniformly in $\alpha$,
  $0\leq \alpha \leq \alpha_M$, and the result follows.  $\ \ \blacksquare$

For the Lipschitzness and Fr\'{e}chet differentiability of the semiflow 
$S_{\alpha}$ (Lemma \ref{lemma28}), 
we first need the estimates of the bilinear operator:
  \begin{lemma}\label{cor3}
    For any $V$ in $H^{s}$, $W$ in $H^{s +1}$, and $s\geq 0$,
    one has
    \begin{itemize}
      \item[(i)\,] $\norm{B_{\alpha}(V,W)}_{s}\leq C(s)
        \norm{V}_{s}\ \norm{W}_{s +1}$.
      \item[(ii)] $\abs{<B_{\alpha}(V,W),\ A^{s}W>}
        \leq D(s) \norm{V}_{s} \norm{W}_{s}^2.$
    \end{itemize}
    Here $C(s)$ and $D(s)$ are constants, which depend on $s$.
 \end{lemma}
 \noindent\textbf{Proof.}
 \begin{itemize}
  \item[(i)] This estimate comes from the inequality (6.6) in \cite{Kim1}. 
     Let, for each fixed wave number $n$,
  \begin{equation}
    (B_{\alpha}(V,W))_n = \sum_{k+m=n} 
                Q_{kmn}(V_k,W_m), \nonumber
  \end{equation}
  where 
  \begin{equation}
  Q_{kmnl}(V_k,W_m) = iP_n\sum_{k+m=n} \left[
         \left( (\mathcal{R}_{\alpha}V)_k\cdot\check{m}\right) V_m +
             V_k^{(j)}\check{m}(\mathcal{R}_{\alpha}W)_m^{(j)} \right]. \nonumber
\end{equation}
 Then 
\begin{eqnarray*}
  \abs{ (\mathcal{R}_{\alpha}V)_k\cdot\check{m} W_m } &\leq&
    \abs{(\mathcal{R}_{\alpha}V)_k}\ \abs{\check{m}}\ \abs{W_m} \\
    &\leq& \frac{1}{1+\alpha^2\abs{\check{k}}^2} \abs{V_k}\ \abs{\check{m}}\ \abs{W_m} \\
    &\leq& \abs{\check{m}}\ \abs{V_k}\ \abs{W_m} \\
  \abs{V_k^{(j)}\check{m}(\mathcal{R}_{\alpha}W)_m^{(j)}}
   &\leq& \abs{V_k^{(j)}}\ \abs{\check{m}}\ \abs{(\mathcal{R}_{\alpha}W)_m^{(j)}} \\
   &\leq& \abs{V_k^{(j)}}\ \abs{\check{m}}\ \frac{1}{1+\alpha^2\abs{\check{m}}^2}
      \abs{W_m^{(j)}} \\
   &\leq& \abs{\check{m}}\ \abs{V_k}\ \abs{W_m}.
\end{eqnarray*}
Since the summation has finitely many terms for each fixed $n$, the bilinear function
has the following inequality,
\[
  \abs{ Q_{kmn}(V_k,W_m) } \leq C \ \abs{\check{m}}\ \abs{V_k}\ \abs{V_m}\ \ 
\]
for an absolute constant $C$. The result follows with $C = C(s)$ in the $H^s$-topology.
  \item[(ii)] This is Lemma 5.3   in \cite{Kim1} with $s=\beta$.\ \ \ $\blacksquare$
 \end{itemize}

  \begin{lemma}\label{lemma28}
    The semiflow $S_{\alpha}(t)v$ is Lipschitz from $[0,\ T]\times B_{\beta}$ into
    $B_{\beta}$ for any $T$ fixed, $T>0$, and it is uniformly Fr\'{e}chet differentiable
    with respect to $v\in B_{\beta}, 0\leq t\leq T$; the above
    properties are uniform in $\alpha$.
  \end{lemma}
  \noindent\textbf{Proof.}
    We closely follow the methodology of Temam \cite{Temam1}, Ch VI, section 8,
    in the context of our semiflow $S_{\alpha}(t)v$ in $B_{\beta},\ \beta > 5/2$.

    Let $V$, $\tilde{V} \in B_{\beta}$ satisfy the equations:
    \begin{eqnarray*}
      \frac{\partial V}{\partial t}+\nu AV+fP_LJP_L\mathcal{R}_{\alpha}V+
        B_{\alpha}(V,V)=F && V(0)=V^0\\
      \frac{\partial \tilde{V}}{\partial t}+\nu A\tilde{V}+
        fP_LJP_L\mathcal{R}_{\alpha}\tilde{V}+
        B_{\alpha}(\tilde{V},\tilde{V})=F &&  \tilde{V}(0)=\tilde{V}^0.
    \end{eqnarray*}
    \begin{itemize}
      \item[(i)\,] First we show a Lipschitz property of the semiflow
        $S_{\alpha}(t):V(0)\rightarrow V(t)$. We set $W(t)=\tilde{V}(t)-V(t)$ and
        $W^0=\tilde{V}^0-V^0$. The difference $W$ satisfies the equation
        {\small 
        \begin{eqnarray}
          \frac{\partial W}{\partial t}+\nu AW+fP_LJP_L\mathcal{R}_{\alpha}W+
           B_{\alpha}(\tilde{V},W)+B_{\alpha}(W,V)=0, && \label{att1}\\
          W(0)=W^0=\tilde{V}^0-V^0.\hspace{1in} &&
        \end{eqnarray}
        }
     Taking $H^{\beta}$-inner product (\ref{att1}) with $W$ and using
     Lemma \ref{cor3} yield
     {\small 
     \begin{equation}\label{att2}
       \frac{1}{2}\frac{d}{dt}\norm{W}_{\beta}^2+\nu\norm{W}_{\beta +1}^2
       \leq c_1\norm{\tilde{V}}_{\beta}\ \norm{W}_{\beta}^2+c_2 \norm{W}_{\beta}^2\
         \norm{V}_{\beta +1}.
     \end{equation}
     }
     By Gronwall's lemma
     \begin{equation}\label{att3}
      \norm{W(t)}_{\beta}^2\leq \norm{W(0)}_{\beta}^2\  e^{\int_0^T2G_{\beta}(\tau )\, d\tau}
     \end{equation}
     where $G_{\beta}(t) =  [ c_1\norm{\tilde{V}}_{\beta}+c_2\norm{V}_{\beta +1}]$.

     This shows the Lipschitz continuity of the semiflow $S_{\alpha}(t)$ with Lipschitz constant
     $C=[\exp (\int_0^T2G_{\beta}(\tau )\,d\tau )]^{1/2}$.
     \item[(ii)] Now we show that the Fr\'{e}chet differentiability of the semiflow 
     $S_{\alpha}(t)$. Consider
       the linearized equations of 3D RNS-$\alpha$ equations
       {\small
       \begin{eqnarray}
         \frac{\partial Z}{\partial t}+\nu AZ +fP_LJP_L Z +
           B_{\alpha}(V,Z)+B_{\alpha}(Z,V)=0, && \label{att4}\\
             Z(0)=Z^0=\tilde{V}^0-V^0.\hspace{1in}
       \end{eqnarray}
       }
      Let $\varphi (t) = \tilde{V}(t)-V(t)-Z(t) = W(t)-Z(t)$. Clearly,
      $\varphi$ satisfies
      {\small
      \begin{equation}\label{att5}
        \frac{\partial \varphi}{\partial t}+\nu A\varphi +fP_LJP_L \varphi +
           B_{\alpha}(V,\varphi)+B_{\alpha}(\varphi,V)
           +B_{\alpha}(W,W)=0, 
      \end{equation}
      }
      with $\varphi(0)=0$.
      This is the exact equations for higher order error $\varphi (t)$.
       Take $H^{\beta}$-inner product (\ref{att5}) with $\varphi$ and use
      Lemma \ref{cor3} to get
      {\small
      \begin{eqnarray}
      \frac{1}{2}\frac{d}{dt}\norm{\varphi}_{\beta}^2+\nu\norm{\varphi}_{\beta +1}^2
      &\leq& c_1\norm{V}_{\beta}\ \norm{\varphi}_{\beta}^2+c_2\norm{V}_{\beta +1}\
           \norm{\varphi}_{\beta}^2 \nonumber\\
           & & \ \ \ + \norm{W}_{\beta}\ \norm{W}_{\beta +1}\
           \norm{\varphi}_{\beta}    \nonumber \\
       &\leq&
       \norm{\varphi}_{\beta}^2 G_{\beta}(t)+\frac{\nu}{2}\norm{\varphi}_{\beta}^2\nonumber \\
      & &  \ \ \  +\frac{1}{2\nu}\norm{W}_{\beta}^2\ \norm{W}_{\beta +1}^2 \label{att6}
      \end{eqnarray}
     }
      Note that, from (\ref{att2}) and (\ref{att3}),
      {\small
      \begin{eqnarray}
        \int_0^T\norm{W(\tau )}_{\beta +1}^2\,d\tau &\leq&
          \frac{1}{\nu}\norm{W^0}_{\beta}^2+\frac{1}{\nu}\norm{W}_{\beta}^2
          \int_0^T 2G_{\beta}(\tau )\,d\tau \nonumber \\
          &\leq&
          \frac{1}{\nu}\norm{W^0}_{\beta}^2
          \left[ 1+ e^{\int_0^T2G_{\beta}(\tau )\,d\tau }\int_0^T 2G_{\beta}(\tau )\,d\tau
          \right] \nonumber \\
          &\equiv& \frac{1}{\nu}\norm{W^0}_{\beta}^2H(T).\nonumber
      \end{eqnarray}
      }
   Substituting this into (\ref{att6}) we obtain
      \begin{equation}
        \norm{\varphi}_{\beta}^2 \leq
          \frac{1}{\nu^2}\norm{W^0}_{\beta}^4\ H(T)
            e^{\int_0^T 4G(\tau )\,d\tau} .
      \end{equation}
      This impiles
      \begin{equation}
        \frac{\norm{\tilde{V}(t)-V(t)-Z(t)}_{\beta}^2}{\norm{\tilde{V}^0-V^0}_{\beta}^2}
        \leq
          \frac{1}{\nu^2}\norm{W^0}_{\beta}^2H(T) e^{\int_0^T 4G(\tau )\,d\tau},
      \end{equation}
      and this shows the Fr\'{e}chet differentiablity of the semiflow
      $S_{\alpha}(t)$. $\ \ \blacksquare$
     \end{itemize}

 \begin{remark} The exponential attractor $\mathcal{M}_{\beta}$ lies in
 the absorbing ball $B_{\beta}$, which is absorbing the initial ball
 $B_{\gamma I}$ in the topology of $H^{\beta}$.
 In that sense, $\mathcal{M}_{\beta}$ is called an
 ``$H^{\beta}$-$H^{\gamma}$" exponential attractor, following
 the usage from damped Hyperbolic PDE's (\cite{Eden}). Technically,
 the global attractor $\mathcal{A}_{\beta}$ is unique in $H^{\beta}$
 in this ``$H^{\beta}$-$H^{\gamma}$" sense.
\end{remark}

 The question arises as to whether $\mathcal{A}_{\beta}$ is the global attractor for
 more general initial data in $H^{\beta}$. The following shows that this is
 true in a sense established in the proof below.
 \begin{corollary}
   The compact global attractor $\mathcal{A}_{\beta}$ attracts
   trajectories with initial data in $H^{\beta},\ \beta > 5/2$,
   with fractal dimension
   uniform in $\alpha$.
 \end{corollary}
 \noindent\textbf{Proof.}
   Let $\epsilon$ be given. Then for every $V_s(0)$ in some arbitrary
   $B_{\gamma I}\subset H^{\gamma},\ \gamma > \beta +4, \ \beta > 5/2$,
   there exists a time $T$ such that:
   \begin{equation}\label{attractor1}
     d_{h,\beta}(V_s(T), \mathcal{A}_{\beta})\leq \epsilon /2,
     \ \mbox{in}\ H^{\beta},
   \end{equation}
   where $d_{h,\beta}(x,y)=\inf_{y\in Y}\norm{x-y}_{\beta}$;
   this follows from $\mathcal{A}_{\beta}$ being a global compact attractor
   in the ``$H^{\beta}$-$H^{\gamma}$" sense.
   We then take any trajectory in some initial ball $\tilde{B}_{\beta}$ in
   $H^{\beta}$ of radius $\tilde{M}_{\beta}$, exactly as in
   Theorem 6.5  and Corollary 6.6  in \cite{Kim1}. For such $V(t)$ with
   $V(0)$ in $\tilde{B}_{\beta}$, we carefully follow the proof of Theorem 6.5
   and Corollary 6.6  in \cite{Kim1}, with $0\leq t\leq T$, $T$ given in (\ref{attractor1}),
   and we can construct $\frac{\eta}{2}=\frac{1}{C_0}\frac{\epsilon}{2}$ to obtain
   \begin{equation}
     \norm{V_s(t)-V(t)}\leq \frac{C_0\eta}{2}=\frac{\epsilon}{2}
   \end{equation}
   on $0\leq t\leq T$. This is actually ``shadowing" of the $V(t)$ trajectory  by
   a $V_s(t)$-trajectory in $H^{\beta}$ on $[0,\ T]$.
   Finally,
   \begin{eqnarray*}
     d_{h,\beta}(V(T),\mathcal{A}_{\beta}) &\leq&
       d_{h,\beta}(V_s(T),\mathcal{A}_{\beta})+\norm{V_s(T)-V(T)}_{\beta} \\
       &\leq& \frac{\epsilon}{2}+\frac{\epsilon}{2} = \epsilon.\ \ \ \blacksquare
   \end{eqnarray*}
   Of course, the rate of attraction is not uniform in the initial data
   as it is to be expected for a topological global attractor; such an attractor
   may attract at a arbitrary slow rate.
\smallskip

\section{Concluding remark}

We constructed exponential attractors in two different functional settings. The exponential
attractor in the Hilbert space (Theorem 3.6 in the Section 3), $\mathcal{M}_\alpha$ lies
in the absorbing ball $B_{\beta}$, $\beta > 5/2$ obtained in the section 2 (Theorem 2.1) and attracts
solution trajectories in $L^2$-topology. In the Banach space, the exponential attractors
$\mathcal{M}_{\beta}$ also lies in the same  $B_{\beta}$ but attracts solution trajectories
in $H^{\beta}$-topology. Two questions arise:

\begin{itemize}
 \item Do they present
     the same or similar asymptotic dynamics for the given system?
 \item Exponential attractors are not unique. Each exponential attractor 
 possesses a unique global attractor, which is a minimal compact attracting sets. 
 Consequently, the intersection of any two exponential attractors is a compact attracting set, 
 for it contains the global attractor. Accordingly, we may ask questions:  
 ``Is the intersection of two exponential attractors still an exponential attractor?".  
 ``How much impact is on the attraction rate?"
   \cite{Pata}
 \end{itemize}



\end{document}